\documentclass[11pt]{amsart}
\usepackage[T1]{fontenc}
\usepackage{lmodern}
\usepackage[a4paper,margin=29mm,headheight=14pt]{geometry}
\usepackage{amsmath,amssymb,amsthm,mathtools}
\usepackage{microtype}
\usepackage{enumitem}
\usepackage[hidelinks]{hyperref}
\usepackage{bookmark}
\usepackage{array}
\numberwithin{equation}{section}
\newtheorem{theorem}{Theorem}[section]
\newtheorem{proposition}[theorem]{Proposition}
\newtheorem{lemma}[theorem]{Lemma}
\newtheorem{corollary}[theorem]{Corollary}
\newtheorem{conjecture}[theorem]{Conjecture}

\newtheorem*{definition*}{Definition}
\theoremstyle{remark}
\newtheorem{remark}[theorem]{Remark}
\newcommand{\R}{\mathbb R}
\newcommand{\dd}{\,d}
\newcommand{\ip}[2]{\langle #1,#2\rangle}

\newcommand{\E}{\mathbb E_{\mu}}
\newcommand{\Q}{\mathcal Q}
\newcommand{\LJac}{L_{\mathrm{Jac}}}
\newcommand{\Id}{\mathrm{Id}}
\DeclareMathOperator{\tr}{tr}
\DeclareMathOperator{\sgn}{sgn}

\setlist[enumerate]{label=\textup{(\roman*)},leftmargin=2em}
\allowdisplaybreaks[1]
\hypersetup{pdftitle={The symmetric Gaussian isoperimetric problem in the plane},
 pdfauthor={Emanuel Milman},pdfkeywords={Gaussian isoperimetry, central symmetry, unconditional sets, curvature, stability}}

\makeatletter
\@namedef{subjclassname@2020}{\textup{2020} Mathematics Subject Classification}
\makeatother

\title[Symmetric Gaussian isoperimetry in the plane]
{The symmetric Gaussian isoperimetric problem in the plane}
\author{Emanuel Milman}
\thanks{Department of Mathematics, Technion - Israel Institute of Technology, Haifa 32000, Israel. Email: emilman@tx.technion.ac.il}
\thanks{The research leading to these results is part of a project that has received funding from the European Research Council (ERC) under the European Union's Horizon 2020 research and innovation programme (grant agreement No 101001677).}

\begin{document}
\begin{abstract}
We resolve the symmetric Gaussian isoperimetric conjecture in the plane:
every minimizer of Gaussian perimeter among centrally symmetric sets of prescribed Gaussian measure is a centered ball, a symmetric strip,
or the complement of one of these sets. A quadrant-reflection argument reduces
the problem to unconditional sets. The planar structural results of
Barchiesi, Cagnetti and Julin then leave only a noncircular strictly convex
stationary oval. We exclude this oval by an explicit volume-preserving and perimeter-decreasing linear
deformation, which geometrically corresponds to an ``unconditional squeezing''. More generally, without requiring any symmetry assumptions, we show that any noncircular closed embedded $\lambda$-self-shrinker is always unstable via linear deformations. We also record
a reduction to the unconditional case in dimension three, using an orthogonal equipartition theorem of Musin.
\end{abstract}

\maketitle

\section{Introduction}

Let $\gamma_n$ denote standard Gaussian measure on $\R^n$, with density
\[
 \rho_n(x)=(2\pi)^{-n/2}\exp(-|x|^2/2).
\]
The \emph{Gaussian volume} of a measurable set $E\subset\R^n$ is
\[
 \gamma_n(E)=\int_E\rho_n(x)\dd x.
\]
For a set of locally finite perimeter, its \emph{Gaussian perimeter} is
\begin{equation}\label{eq:perimeter-def}
 P_{\gamma_n}(E)=\int_{\partial^*E}\rho_n(x)\dd\mathcal H^{n-1}(x),
\end{equation}
where $\partial^*E$ denotes its reduced boundary (of course $\partial^*E = \partial E$ for a set with smooth boundary). For a measurable set which is not of
locally finite perimeter, we set $P_{\gamma_n}(E)=+\infty$.
All sets and symmetries are understood modulo null sets. We refer to
\cite{AFP,Maggi} for the underlying facts about sets of finite perimeter. We denote by $B_r^k$ the
Euclidean ball of radius $r$ centered at the origin in $\R^k$. 

The classical Gaussian isoperimetric theorem of Sudakov--Tsirel'son
\cite{SudakovTsirelson} and Borell \cite{Borell} states that half-spaces minimize
Gaussian perimeter given prescribed Gaussian volume. In this work we impose
the additional constraint of \emph{central symmetry}, $E=-E$. A second
more restrictive symmetry class which will be useful for us is that of \emph{unconditional} sets, namely sets invariant under reflection across each coordinate hyperplane.

The question of identifying the Gaussian isoperimetric minimizers among all centrally symmetric sets was put forth by Franck Barthe in \cite{Barthe}. 
For unconditional sets and $\ell_\infty^n$ anisotropic Gaussian perimeter, Barthe resolved the corresponding isoperimetric problem: the minimizers are precisely slabs or their complements. However, the problem for the usual isotropic Gaussian perimeter is more complicated even in the planar case, where the centered Euclidean ball can have smaller Gaussian perimeter than the corresponding slab (or its complement) with the same Gaussian volume. The natural conjecture is the following one (first explicitly appearing in the work by Heilman \cite[Conjecture 1.3]{Heilman}). Chakrabarti--Regev \cite{ChakrabartiRegev} and O'Donnell \cite{ODonnell} have also proposed variants of the problem for Gaussian noise stability.

\begin{conjecture}[Symmetric Gaussian isoperimetric conjecture]
\label{conj:Barthe}
Let $n \geq 1$. For any $v\in(0,1)$, the minimum of $P_{\gamma_n}(E)$ over
centrally symmetric sets with $\gamma_n(E)=v$ is attained by a round cylinder
\[
 B_r^k\times\R^{n-k}
\]
or its complement, for an appropriate $r>0$ and $1 \leq k \leq n$. Moreover, up to rotations (and null set modifications), these are the unique minimizers. 
\end{conjecture}

Note that the rank $k$ is allowed to depend on the Gaussian volume $v$ and the ambient dimension $n$. 
Of course, complementation, null set modifications and orthogonal transformations all preserve the perimeter.
In dimension $n=1$ the only minimizing configurations are centered intervals
and pairs of symmetric tails: the interval minimizes above volume $1/2$,
the tails below volume $1/2$, and both at volume $1/2$; see, for example,
\cite[Appendix~A]{BCJ}.

\medskip

Several partial and related results for centrally symmetric sets have been obtained in the literature.
Lata{\l}a and Oleszkiewicz \cite[Theorem~3]{LatalaOleszkiewicz} proved that
symmetric strips minimize Gaussian perimeter multiplied by the centered
inradius among centrally symmetric convex sets of fixed Gaussian volume.
This is a width-weighted inequality, rather than the unweighted conjecture.
Barthe's aforementioned result \cite{Barthe} concerns an anisotropic
perimeter, rather than \eqref{eq:perimeter-def}. For the isotropic Gaussian
perimeter, La Manna \cite[Theorem~1 and Proposition~1]{LaManna} proved,
for $n\ge2$,
quantitative $L^1(\gamma_n)$-local minimality of centered balls of radius
$0<r<\sqrt{n+1}$, and instability for $r>\sqrt{n+1}$. 
Barchiesi and Julin \cite[Corollary~1]{BarchiesiJulin} proved that symmetric
strips are the unique minimizers, up to rotations and null sets, for volumes
sufficiently close to one, with a dimension-independent threshold.
Complementation gives the corresponding result near zero. 
Heilman \cite{Heilman} obtained partial classification results in arbitrary dimension: for a minimizer $E$ such that $E$ or its complement is convex, he proved
that its boundary is a centered round cylinder under appropriate average curvature conditions. 
A distinct but related problem involving a Gaussian-weighted total-curvature functional was studied by
Fusco and La Manna \cite{FuscoLaManna}.

\medskip

These results do not cover all curvature or volume regimes, and leave even the planar
problem unresolved. The difficulty is due to the fact that stationarity alone does not characterize the conjectured minimizers even in the plane: noncircular embedded centrally symmetric Gaussian stationary curves were constructed by Chang \cite{Chang}. An interesting attempt to reduce the problem to the planar case was recently suggested by Barchiesi, Cagnetti and Julin \cite{BCJ}: they introduced a
double-arc rearrangement \emph{for unconditional sets}, proved that it reduces Gaussian perimeter for such sets, and obtained a useful planar
structural result which we shall use in this work.

\bigskip

Our main result gives the full planar classification.

\begin{theorem}[Confirmation in the plane]
\label{thm:main}
Conjecture \ref{conj:Barthe} holds for $n=2$. 
\end{theorem}

The proof has two separate ingredients. First, reflection of a suitable
equipartitioned quadrant gives an unconditional competitor of no greater
perimeter. The reflection argument in the planar case is elementary, but surprisingly also works in dimension $n=3$ (see Section \ref{sec:unconditional}). 
Combined with the planar results of \cite{BCJ}, this reduces the
problem to disqualifying a strictly convex stationary oval whose curvature has primitive
period $\pi$ in the normal angle. 

Second, we prove that this oval is unstable; more precisely, that it has a
negative second variation of perimeter along a volume-preserving diagonal linear
deformation (which geometrically corresponds to an ``unconditional squeezing''). In fact, we develop a general instability criterion which does not require any symmetry assumptions. The decisive analytical estimate is the Key Inequality of Proposition
\ref{prop:key} below. Its proof uses only the curvature ODE, the sharp Poincar\'e (or Wirtinger) inequality on $[0,\pi]$, Cauchy--Schwarz, and an exact integration-by-parts identity.

\medskip

Our results also admit an interpretation and extension in terms of general
$\lambda$-self-shrinkers, defined by the equation
\[
 k-\langle x,N\rangle \equiv \lambda,
\]
where $N$ is the outward unit normal and $k$ is the signed curvature,
positive on convex curves. This is precisely the stationarity equation
for Gaussian perimeter under a Gaussian-volume constraint; the case
$\lambda=0$ corresponds, after a fixed rescaling, to the usual
self-shrinkers for curve-shortening flow \cite{ChengWei}. Abresch and Langer \cite{AbreschLanger} showed that centered circles are the only closed embedded self-shrinkers (with $\lambda = 0$). However, as already mentioned, Chang \cite{Chang} showed that noncircular closed embedded
$\lambda$-self-shrinkers exist for certain $\lambda < 0$. 
Moreover, he showed that when $a = -\lambda \in (0,2 / \sqrt{3})$, there exist centrally symmetric solutions (and an inspection of his proof verifies that in fact their \emph{primitive} period is $\pi$, meaning that they are not invariant under $2\pi/m$-rotations for $m\geq 3$). We observe that this is in fact an equivalence, and show that every noncircular closed embedded solution is unstable as follows:

\begin{theorem} \label{thm:main2}
Let $\Sigma$ be a noncircular closed embedded $\lambda$-self-shrinker in $\R^2$ (no symmetry is assumed). Then:
\begin{enumerate}
\item $\Sigma$ has primitive period $\pi$ if and only if $a = -\lambda \in (0,2 / \sqrt{3})$. 
\item $\Sigma$ is Gaussian unstable via linear deformations.  
If $\Sigma$ is centrally symmetric then the instability is in the centrally symmetric class. 
If $\Sigma$ is unconditional then the instability is in the unconditional class (but possibly non-linear if the primitive period is $< \pi$).
\end{enumerate}
\end{theorem}

Note that without symmetry restrictions, instability already follows from
McGonagle and Ross \cite[Corollary~4.8]{McGonagleRoss}, who used translations to prove that hyperplanes are the only Gaussian stable complete two-sided properly immersed Gaussian stationary hypersurfaces of finite Gaussian area. In contrast, our argument establishes instability via linear
deformations, and also yields the symmetry-preserving conclusions specified above.

\bigskip

The rest of this work is organized as follows. 
Section~\ref{sec:unconditional} proves the unconditional reductions in
dimensions two and three; only the two-dimensional lemma is used to prove
Theorem~\ref{thm:main}. Section~\ref{sec:structure} identifies the remaining
oval. Section~\ref{sec:stability} reduces its instability to the Key
Inequality, which is proved in Section~\ref{sec:key-proof}.
In Section~\ref{sec:completion} our two main theorems are established, including the transfer
of equality cases back to arbitrary centrally symmetric sets.
The appendix provides the exact volume-preserving deformation and the
second-variation computation.

\bigskip

\textbf{AI Declaration.} The reduction of the planar problem to the problem of excluding noncircular stationary ovals by proving their instability
is well known to experts. 
While the reduction to the unconditional case in the planar setting is trivial,  AI was able to extend the argument to the three dimensional case using Musin's theorem. The author asked AI to show that an unconditional oval must be unstable, and surprisingly it failed, attempting instead to rely on numerical arguments. The author then suggested to use linear deformations, which resulted in the emergence of the sufficient condition $S < 0$ for instability. AI could not establish this sufficient condition, and so the author suggested that it should use Poincar\'e's inequality, after which AI was able to establish the Key Inequality and thereby conclude the proof. Essentially the same argument yields Theorem \ref{thm:main2}. AI was also used to draft a preliminary version of this paper, which was revised and refined by the author, who takes full responsibility for everything written.

\section{Reduction to unconditional sets in dimensions two and three}
\label{sec:unconditional}

We first record a standard reflection observation valid in every dimension. For an
open set $U$, write
\[
 P_{\gamma_n}(E;U)=
 \int_{\partial^*E\cap U}\rho_n\dd\mathcal H^{n-1}.
\]

\begin{lemma}[Reflection of an equipartitioned orthant]
\label{lem:reflection}
Let $E\subset\R^n$ have finite Gaussian perimeter. Suppose an orthonormal
basis $\{e_i\}_{i=1,\ldots,n}$ determines open orthants $Q_\varepsilon = \{x \in \R^n : \ip{x}{e_i} \varepsilon_i > 0 \;\; \forall i =1, \ldots, n \}$, $\varepsilon\in\{-1,1\}^n$,
such that
\[
 \gamma_n(E\cap Q_\varepsilon)=2^{-n}\gamma_n(E)
 \quad \;\; \forall \varepsilon \in\{-1,1\}^n . 
\]
Then there exists a measurable $F \subset \R^n$, unconditional in this basis, with
$\gamma_n(F)=\gamma_n(E)$ and $P_{\gamma_n}(F)\le P_{\gamma_n}(E)$.
More precisely, if $F$ denotes the reflection of any $E\cap Q_{\varepsilon}$ into all orthants, then 
\begin{equation}\label{eq:reflection-exact}
 \gamma_n(F)=2^n\gamma_n(E\cap Q_\varepsilon),\qquad
 P_{\gamma_n}(F)=2^nP_{\gamma_n}(E;Q_\varepsilon).
\end{equation}
\end{lemma}

\begin{proof}
Gaussian measure is invariant under orthogonal reflections, proving the
volume identity in (\ref{eq:reflection-exact}). Across each reflecting hyperplane, the one-sided BV
traces of the reflected characteristic functions agree. The BV gluing
formula therefore gives no perimeter contribution on that hyperplane.
Intersections of two or more cutting hyperplanes have codimension at least
two and carry no perimeter measure. This proves the perimeter identity in (\ref{eq:reflection-exact}).
These assertions can be applied locally in balls, where finite Gaussian
perimeter implies finite Euclidean perimeter, and then exhausted over
$\R^n$; see \cite[Chapter~3]{AFP} for the trace and gluing formulas.
Choosing $Q$ of least relative perimeter yields
\[
 P_{\gamma_n}(F)
 =2^nP_{\gamma_n}(E;Q)
 \le\sum_\varepsilon P_{\gamma_n}(E;Q_\varepsilon)
 \le P_{\gamma_n}(E).
\]
Any perimeter of the original set on a cutting hyperplane only strengthens
the final inequality.
\end{proof}

\begin{lemma}[Unconditional reduction in dimension two]
\label{lem:unconditional-2}
For every centrally symmetric set $E\subset\R^2$ of finite Gaussian
perimeter, there is an unconditional set $F\subset\R^2$ satisfying
\[
 \gamma_2(F)=\gamma_2(E),\qquad
 P_{\gamma_2}(F)\le P_{\gamma_2}(E).
\]
\end{lemma}

\begin{proof}
Put $\dd\nu=\mathbf1_E\dd\gamma_2$ and $v=\nu(\R^2)$. This finite
measure is even and vanishes on lines. Let $Q_\theta$ be the open quadrant
between angles $\theta$ and $\theta+\pi/2$, and set
$m(\theta)=\nu(Q_\theta)$. Dominated convergence gives continuity of $m$.
Every half-plane through the origin has $\nu$-measure $v/2$, so
\[
 m(\theta)+m(\theta+\pi/2)=v/2.
\]
Thus $m(\theta)-v/4$ changes sign under a quarter-turn and has a zero.
Evenness gives equality of opposite quadrant masses; at such a zero all
four quadrant masses equal $v/4$. Apply Lemma~\ref{lem:reflection} and
rotate the resulting set to the standard coordinate axes.
\end{proof}

\begin{lemma}[Unconditional reduction in dimension three]
\label{lem:unconditional-3}
For every centrally symmetric set $E\subset\R^3$ of finite Gaussian
perimeter, there is an unconditional set $F\subset\R^3$ satisfying
\[
 \gamma_3(F)=\gamma_3(E),\qquad
 P_{\gamma_3}(F)\le P_{\gamma_3}(E).
\]
\end{lemma}

\begin{proof}
We first introduce some notation. 
The hyperoctahedral group $B_k = A_k \rtimes S_k$ is the
semidirect product of the sign-change group
$A_k = (\mathbb{Z}/2)^k = \langle\lambda_1,\ldots,\lambda_k\rangle$
and the symmetric group $S_k$, where $S_k$ acts by
$\sigma(\lambda_i)=\lambda_{\sigma(i)}$.
The generator $\lambda_i \in A_k$ acts on $\mathbb{R}^k$ by
negating the $i$-th coordinate.
Let $O(k)$ be the orthogonal group, which we identify with the space of ordered orthonormal frames in $\R^k$. The group $B_k$ acts on this space by sign changes and permutations of the frame vectors. 

For each pair $\{i,j\}\subset[k]$ with $i<j$, let
$\mathbb{R}_{\lambda_i+\lambda_j}$ be the one-dimensional real
$A_k$-representation on which $\lambda_\ell$ acts by $-1$ if
$\ell\in\{i,j\}$ and by $+1$ otherwise.  Set
\[
 W_2^+ = \bigoplus_{1\le i<j\le k}\mathbb{R}_{\lambda_i+\lambda_j},
\]
extended to a $B_k$-representation by
$\sigma\,e_{\{i,j\}}=e_{\{\sigma(i),\sigma(j)\}}$.

The key idea is to use Musin's equivariant zero theorem \cite[Theorem~1.1]{Musin} for $k=3$:
every continuous $B_k$-equivariant map from $O(k)$ into $W_2^+$ 
has a zero.

Set $\dd\nu=\mathbf1_E\dd\gamma_3$. For an orthonormal frame
$(e_1,e_2,e_3)$, consider the three coordinates
\begin{equation}\label{eq:pair-sign-moments}
 \int_{\R^3}
 \sgn\ip{x}{e_i}\,\sgn\ip{x}{e_j}\dd\nu(x),
 \qquad 1\le i<j\le3.
\end{equation}
They define a map from $O(3)$ to $W_2^+$ which is clearly $B_3$-equivariant. It is continuous
by dominated convergence, since $\nu$ vanishes on every plane through the
origin. Hence a frame $(e_1,e_2,e_3)$ exists for which all three quantities vanish.

Write $\sigma_i(x)=\sgn\ip{x}{e_i}$ and $v=\nu(\R^3)$. Evenness gives
\[
 \int\sigma_i\dd\nu=0,\qquad
 \int\sigma_1\sigma_2\sigma_3\dd\nu=0.
\]
Outside the cutting planes, the orthant
$Q_\varepsilon=\{x:\varepsilon_i\ip{x}{e_i}>0,\ i=1,2,3\}$ has indicator
\[
 \mathbf1_{Q_\varepsilon}
 =\frac18\prod_{i=1}^3(1+\varepsilon_i\sigma_i).
\]
Integrating its expansion, the odd-degree terms vanish by evenness and
the degree-two terms by \eqref{eq:pair-sign-moments}. Thus every orthant has
mass $v/8$. Lemma~\ref{lem:reflection}, followed by a rotation, completes
the proof.
\end{proof}

\begin{remark}\label{rem:scope-reduction}
For $n=2,3$, these lemmas prove equality of the centrally symmetric and
unconditional isoperimetric infima at every volume. They do not by
themselves imply that every centrally symmetric minimizer is unconditional
up to rotation. The planar equality cases will be treated separately in
Section~\ref{sec:completion}. The reflection operation need not preserve
convexity. Nor does the argument establish a reduction in all dimensions:
when $n\ge4$, vanishing of pairwise sign moments leaves higher even-degree
sign moments uncontrolled.
\end{remark}

\section{The exceptional stationary oval}
\label{sec:structure}

From now on we exclusively restrict to the planar case $n=2$, writing $P_\gamma=P_{\gamma_2}$ and
$\gamma=\gamma_2$. The unconditional problem has a minimizer at every
volume $v\in(0,1)$ by a standard argument. Indeed, a minimizing sequence has uniformly bounded
Euclidean perimeter on each bounded ball $B_R$. BV compactness and a diagonal
argument give local $L^1$ convergence of characteristic functions.
The uniform tail estimate by $\gamma(\R^2\setminus B_R)$ upgrades this
to global $L^1(\gamma)$ convergence. Volume and unconditionality are
preserved, and standard lower semicontinuity of perimeter yields a minimizer (see e.g.~\cite{BCJ}). 

\medskip

We will use the following planar result of Barchiesi, Cagnetti and Julin \cite{BCJ}. 

\begin{proposition}[Planar unconditional minimizers \cite{BCJ}]\label{prop:BCJ}
Let $E\subset\R^2$ minimize Gaussian perimeter among unconditional sets
of prescribed Gaussian volume. After replacing $E$ by its complement
if necessary, and choosing its open representative, one of the following
holds:
\begin{enumerate}
\item $E$ is a symmetric strip;
\item $E$ is a disk centered at the origin;
\item $E$ is bounded and strictly convex, $\partial E$ is a real analytic
simple closed curve with positive curvature $k$ and exterior unit normal $N$, and
\begin{equation}\label{eq:nonradial}
 N(x)\ne x/|x|
 \qquad\text{for every }x\in\partial E\cap(0,\infty)^2.
\end{equation}
\end{enumerate}
\end{proposition}
\begin{proof}
This follows from \cite[Proposition~9]{BCJ} and its proof. Their Proposition~4 gives real analyticity of the entire boundary
in dimension two, and Proposition~6 reduces the disconnected case
to a symmetric strip or its complement. In the connected case,
the proof of Proposition~9 establishes compactness and then,
after possible complementation, strict positivity of curvature (the argument will be repeated below). 
Thus the boundary is a
real-analytic simple closed curve enclosing a strictly convex
domain. Finally, their Corollary~2 shows that a radial outward
normal at any boundary point in the open positive quadrant
forces the domain to be a centered disk.
\end{proof}

Let us describe alternative (iii) more precisely. Let $K$ be its
bounded convex body and $\Sigma=\partial K$. A standard argument  
using the first variation of Gaussian perimeter
for the fixed Gaussian volume constraint (see e.g.~\cite{EMilmanNeemanAnnals,EMilmanCubeCPAM}) verifies that the Gaussian-weighted mean-curvature of $\Sigma$ must be constant:
\begin{equation}\label{eq:stationarity}
 \exists \lambda \in \R  \qquad k-\ip{x}{N} \equiv \lambda\qquad\text{on }\Sigma.
\end{equation}
This is the \emph{Gaussian stationarity equation}. Although the minimization
is restricted to unconditional sets, the equation holds pointwise:
averaging a test function over the (discrete) coordinate-reflection group does not
change its pairing with the invariant Euler--Lagrange expression. 

\begin{definition*}[$\lambda$-self-shrinkers and stationary ovals]
\hfill
\begin{itemize}
\item A smooth embedded (connected) curve $\Sigma$ satisfying (\ref{eq:stationarity}) is called a $\lambda$-self-shrinker. 
\item When $\Sigma$ is in addition closed (i.e.~compact) and noncircular, it is called a stationary oval. 
\end{itemize}
\end{definition*}
\noindent Note that there is no assumption of central symmetry in either definition, to accommodate the general formulation of Theorem \ref{thm:main2}.  

\medskip

As alluded to in the proof of Proposition \ref{prop:BCJ}, a stationary oval is always necessarily strictly convex.  Indeed, differentiating (\ref{eq:stationarity}) along arclength $s$ (counterclockwise), with $T$ denoting the positively oriented unit tangent vector, yields
\begin{equation} \label{eq:partial-s}
\partial_s x=T,\quad\partial_s N=k T , \quad \partial_s T = - k N , \quad \partial_s k=k\ip{x}{T} ,
\end{equation}
and the latter integrates to 
 \begin{equation} \label{eq:k-rho}
 k \exp(-|x|^2/2) \equiv C \text{ on } \Sigma . 
 \end{equation}
Using the turning-number theorem $\int_{\Sigma} k \, ds = 2 \pi$, it follows that $C > 0$ and hence $k > 0$ as well. 

Being strictly convex, we parametrize the curve $[0,2\pi) \ni \theta \mapsto x(\theta) \in \Sigma$ by its outward normal angle $\theta$, and denote the outward unit normal vector, unit tangent vector, and support function (respectively) by
\[
 N(\theta)=(\cos\theta,\sin\theta),\qquad
 T(\theta)=(-\sin\theta,\cos\theta),\qquad
 h(\theta)=\ip{x(\theta)}{N(\theta)}.
\]
Primes will henceforth always denote differentiation with respect to $\theta$.
Some standard useful identities are
\begin{equation}\label{eq:support}
 x=hN+h'T,\qquad
  x'=(h+h'') T = \frac1k T ,\qquad
 \dd s=\frac{\dd\theta}{k} .
\end{equation}
With $a=-\lambda$, equation \eqref{eq:stationarity} says $h=k+a$, and consequently $h + h'' = \frac1k$ translates to
\begin{equation}\label{eq:curvature-ode}
 k''+k+a=\frac1k.
\end{equation}

\begin{remark}We will see in Section \ref{sec:key-proof} that the ODE (\ref{eq:curvature-ode}) is explicitly integrable. \end{remark}

A centrally symmetric stationary oval clearly has period $\pi$, and as the origin belongs to the interior of $K$, it follows that $h > 0$.

\subsection{Primitive period \texorpdfstring{$\pi$}{pi}}

We conclude this section with a lemma which is not essential for the proof of Theorem \ref{thm:main}, but nevertheless is 
the most convenient option. 
It states that not only does the remaining oval in Proposition \ref{prop:BCJ} have period $\pi$, but in fact $\pi$ is its minimal (or primitive) period.

\begin{lemma}[Primitive period of the remaining oval]\label{lem:period}
In alternative \textup{(iii)} of Proposition~\ref{prop:BCJ}, the curvature
$k$ is a positive nonconstant solution of \eqref{eq:curvature-ode} with
primitive period exactly $\pi$. After interchanging the coordinate axes if
necessary, it satisfies
\begin{equation}\label{eq:phase}
 k'(0)=k'(\pi/2)=0,
 \qquad k'<0\quad\text{on }(0,\pi/2).
\end{equation}
\end{lemma}

\begin{proof}
Central symmetry gives $h(\theta+\pi)=h(\theta)$, and the coordinate
reflections give $h(-\theta)=h(\theta)$ and
$h(\pi-\theta)=h(\theta)$. These statements also hold for $k=h-a$;
in particular $k'(0)=k'(\pi/2)=0$.

For $0<\theta<\pi/2$, the support point $x(\theta)$ belongs to the open
positive quadrant. Indeed, if $N_i>0$ and $x_i<0$, reflection of the
support point in the $i$th coordinate would increase its scalar product
with $N$. If $N_i>0$ and $x_i=0$, reflection would give two different
outer normals at the same point, contradicting smoothness. By
\eqref{eq:support} and $h>0$, the equality $N=x/|x|$ is equivalent to
$h'=0$. Thus \eqref{eq:nonradial} says precisely that $k'$ never vanishes
on $(0,\pi/2)$. 

Reflection about $\pi/2$ gives
$k'(\pi-\theta)=-k'(\theta)$, so $k'$ also never vanishes on
$(\pi/2,\pi)$. Its only zeros in $[0,\pi)$ are consequently $0$ and
$\pi/2$. If $\tau \in(0,\pi)$ were a period of $k$, then
$k'(\tau)=k'(0)=0$, forcing $\tau=\pi/2$. But strict monotonicity on
$(0,\pi/2)$ gives $k(\pi/2)\ne k(0)$, a contradiction. Since $\pi$
is already a period by central symmetry, it is the primitive period.
Finally, interchanging the coordinate axes replaces $k(\theta)$ by
$k(\pi/2-\theta)$ and reverses the sign of its derivative on
$(0,\pi/2)$, giving \eqref{eq:phase}.
\end{proof}

\begin{remark}[An alternative nodal-domain route]\label{rem:Heilman}
The infinitesimal-rotation argument in \cite[\S11.2, Corollary~11.9]{Heilman} offers another way to eliminate centrally-symmetric stationary ovals with higher frequency without invoking unconditionality. The idea is as follows: on such ovals, the rotational Jacobi field is $\ip{Jx}{N}=-h'=-k'$, where $J$ is a quarter-turn.
If the oval were stable and had more than four nodal arcs, one could construct a nonzero
even function $g$ by taking a linear combination of two
antipodal pairs of nodal truncations, with coefficients chosen so that
$g$ has zero Gaussian boundary mean. This function belongs to $H^1$ and vanishes on an unused open arc. 
It follows that $Q(g)=0$ and hence by stability that $L_{Jac} g \equiv C = 0$ weakly (see the definition of $Q$ and $L_{Jac}$ in Section \ref{sec:stability}). Uniqueness for the resulting second-order linear ODE yields $g\equiv0$, a contradiction. This leaves at most four nodal arcs, hence primitive period $\pi$. 
To obtain a complete structural reduction via this route, one would
also need to justify regularity at the origin (where the central
symmetry has a fixed point), handle possible disconnected boundaries
(including pairs of components exchanged by central symmetry), and
establish that each remaining component is a line or a compact
strictly convex curve.
We instead prefer to use Proposition~\ref{prop:BCJ}, combined with the
elementary reduction to unconditional sets in
Lemma~\ref{lem:unconditional-2}.
\end{remark}

\section{The Jacobi operator and the linear instability test}
\label{sec:stability}

Let $\Sigma = \partial K$ be a stationary oval as in the previous section. It satisfies the Gaussian stationarity equation (\ref{eq:stationarity}), or equivalently, 
recalling that $a = - \lambda$, the curvature ODE (\ref{eq:curvature-ode}). Unless otherwise stated, we do not assume here that $K$ is an isoperimetric minimizer or even centrally symmetric, and develop a general criterion for testing the instability of $\Sigma$, which will be used in the proofs of Theorems \ref{thm:main} and \ref{thm:main2}. 

\subsection{Second variation of perimeter and the Jacobi operator}

A standard idea in the calculus of variations is to construct competitors to $\Sigma$ by considering variations with normal velocity $f$. 
We refer to  e.g.~\cite{EMilmanNeemanAnnals,EMilmanCubeCPAM} for the standard facts and formulas used below. For completeness, we also outline the computation in the appendix, as it is particularly simple in our smooth planar setting.

Define the normalized Gaussian probability measure on $\Sigma$
\[
 \dd\mu=\frac{\rho_2\dd s}{\int_{\Sigma} \rho_2 \dd s} .
\]
Recalling that $\rho_2$ is proportional to $1/k$ by (\ref{eq:k-rho}), and that 
$\dd s=k^{-1}\dd\theta$, it follows that for every function $f$ on $\Sigma$ we have
\begin{equation}\label{eq:mu}
 \E f=
 \frac{\displaystyle\int_0^{2\pi} f(\theta) k(\theta)^{-2} d\theta}
 {\displaystyle\int_0^{2\pi} k(\theta)^{-2} d\theta} .
\end{equation}
Clearly a normal velocity
$f$ has zero first Gaussian-volume variation if and only if $\E f=0$. When $K$ is centrally symmetric, the evenness of $f$ corresponds to $\pi$-periodicity in $\theta$-coordinates, and both integrals in (\ref{eq:mu}) may then be taken over $[0,\pi]$. 

The stationarity equation (\ref{eq:stationarity}) precisely characterizes the property that $\Sigma$ is a critical (or stationary) point for Gaussian perimeter for any deformation fixing the Gaussian volume constraint, or equivalently, that the first variation of Gaussian perimeter is a constant multiple $\lambda$ (the Lagrange multiplier) of the  first variation of Gaussian volume. As for the second variation of Gaussian perimeter for a fixed Gaussian volume constraint, it is governed by the Gaussian Jacobi operator, defined as
\begin{equation}\label{eq:Jacobi-intrinsic}
 \LJac f := \Delta_{\Sigma,\gamma} f +(1+k^2)f ,
\end{equation}
where
\[
\Delta_{\Sigma,\gamma}:=\partial_s^2-\ip{x}{T}\partial_s
\]
 is the Gaussian-weighted Laplacian on $\Sigma$ (of course when $\Sigma$ is just a curve, there is no need to sum over partial derivatives). 
Since $\partial_s=k\partial_\theta$ and $\ip{x}{T}=h'=k'$, the first-order
terms cancel, giving the following slick representation in $\theta$-coordinates:
\begin{equation}\label{eq:Jacobi-angle}
 \LJac f=k^2f''+(1+k^2)f.
\end{equation}
For smooth $f,g$, the corresponding (normalized) second-variation (or index) form is
\begin{equation}\label{eq:Q}
 \begin{aligned}
 \Q(f,g) :=-\E[f\LJac g] =\frac{\displaystyle\int_0^{2\pi}
 \bigl(f'g'-(1+k^{-2})fg\bigr)\dd\theta}
 {\displaystyle\int_0^{2\pi} k^{-2}\dd\theta},
 \qquad \Q(f) :=\Q(f,f).
 \end{aligned}
\end{equation}
The last transition is immediate from integration by parts with periodic boundary conditions, establishing the symmetry of $Q$ (and corresponding self-adjointness of $L_{Jac}$ on $L^2(\mu)$). 

Now every smooth $f$ with $\E f=0$ is the initial normal velocity
of an exactly Gaussian-volume-preserving
deformation, along which the second derivative of perimeter is
$P_\gamma(K)\Q(f)$; see the appendix. 
When $K$ and $f$ are centrally symmetric or unconditional, this deformation preserves the same symmetry.
 Consequently, if $\Sigma$ minimizes Gaussian perimeter under a Gaussian volume constraint, whether without symmetry restrictions or in the class of centrally symmetric or unconditional sets, it must be stable in the following sense:
\begin{equation} \label{eq:stability}
\E f =0 \;\; \Rightarrow \;\; \Q(f) \geq 0 ,
\end{equation}
for all smooth, even smooth or unconditional smooth functions $f$, respectively. We refer to the appendix for more information. 

\subsection{Linear deformations}

For a constant symmetric $2\times2$ matrix $B$, set
\[
 f_B=\ip{Bx}{N},
\]
and note that $h=f_{\Id}$. This is the normal velocity of the linear deformation $e^{tB} K$ generated by the map $x \mapsto e^{tB} x$. 
If $K$ is centrally symmetric, $f_B$ is even; if $K$ is unconditional and $B$ is diagonal, $f_B$ is
unconditional.

To compute the resulting second variation, first compute using \eqref{eq:stationarity} and (\ref{eq:partial-s}) that
\[
 \Delta_{\Sigma,\gamma}x=-x-\lambda N,
 \qquad \Delta_{\Sigma,\gamma}N=-k^2N ,
 \]
and then apply the product rule for the Laplacian, yielding
\begin{equation}\label{eq:linear-Jacobi}
 \LJac f_B=2k\ip{BT}{T}-\lambda\ip{BN}{N}.
\end{equation}
 
Now choose diagonal and shear matrices
\[
 B_D :=\begin{pmatrix}1&0\\0&-1\end{pmatrix},\qquad
 B_{\mathrm{sh}} :=\begin{pmatrix}0&1\\1&0\end{pmatrix}.
\]
Recalling (\ref{eq:support}), their corresponding normal velocities are computed to be
\begin{equation}\label{eq:linear-speeds}
 f_D=h\cos2\theta-k'\sin2\theta,
 \qquad
 f_{\mathrm{sh}}=h\sin2\theta+k'\cos2\theta.
\end{equation}
Equation \eqref{eq:linear-Jacobi} gives
\begin{equation}\label{eq:linear-images}
 \LJac f_D=(a-2k)\cos2\theta,
 \quad
 \LJac f_{\mathrm{sh}}=(a-2k)\sin2\theta,
 \quad
 \LJac h=2k+a.
\end{equation}

\medskip

For the instability criterion below we assume that $a>0$; we will verify
this for every stationary oval in
Lemma~\ref{lem:period-parameter}. In particular, $h=k+a>0$, without
requiring any symmetry assumption. Since neither linear velocity need
have zero Gaussian boundary mean, we balance both by a dilation:
\begin{equation}\label{eq:balance}
 \begin{aligned}
 c_D&:=\frac{\E f_D}{\E h},
 &g_D&:=f_D-c_Dh,\\
 c_{\mathrm{sh}}&:=\frac{\E f_{\mathrm{sh}}}{\E h},
 &g_{\mathrm{sh}}&:=f_{\mathrm{sh}}-c_{\mathrm{sh}}h.
 \end{aligned}
\end{equation}
Since $h>0$, $c_D, c_{\mathrm{sh}}\in \R$ are well defined, and $\E g_D=\E g_{\mathrm{sh}}=0$. 
The quantity controlling the instability is
\begin{equation}\label{eq:S}
 S:=2\E k^2+a\E k-a^2.
\end{equation}

\begin{lemma}[Second variation identity]\label{lem:balance}
Assume that $a>0$. Then
\begin{equation}\label{eq:balance-identity}
 \Q(g_D)+\Q(g_{\mathrm{sh}})
 =\bigl(1-c_D^2-c_{\mathrm{sh}}^2\bigr)S,
 \qquad c_D^2+c_{\mathrm{sh}}^2<1.
\end{equation}
\end{lemma}

\begin{proof}
For every trace-free symmetric matrix $B$, we have
\begin{equation}\label{eq:Euclidean-cancellation}
 \E[kf_B]=0.
\end{equation}
Indeed, \eqref{eq:k-rho} makes this Gaussian integral a positive constant multiple
of the Euclidean boundary integral of $f_B$, and so by the divergence theorem
\[
 \int_\Sigma f_B\dd s
 =\int_\Sigma\ip{Bx}{N}\dd s
 =\int_K\tr B\dd x=0.
\]
Next, (\ref{eq:linear-speeds}) and (\ref{eq:linear-images}) yield
\begin{align}
 \Q(f_D)+\Q(f_{\mathrm{sh}}) & = - (\E(f_D \LJac f_D) + \E(f_{\mathrm{sh}} \LJac f_{\mathrm{sh}})) \notag \\
 &=-\E\bigl[(a-2k)(f_D\cos2\theta+f_{\mathrm{sh}}\sin2\theta)\bigr]
 \notag\\
 &=-\E[(a-2k)(k+a)]=S.\label{eq:trace-identity}
\end{align}
The same formulas, together with \eqref{eq:Euclidean-cancellation}, yield
\begin{align*}
  \Q(h) &=-\E[(2k+a)(k+a)] = -S - 2 a (\E k + a) = -S - 2 a \E h , \\
 \Q(f_D,h) &=-\E[(2k+a) f_D]=-a\E f_D = - a c_D \, \E h , \\
 \Q(f_{\mathrm{sh}},h) &= -\E[(2k+a) f_{\mathrm{sh}}] = -a\E f_{\mathrm{sh}} = - a c_{\mathrm{sh}} \, \E h . 
\end{align*}
We conclude that for $* \in \{D , \mathrm{sh} \}$, 
\begin{equation} \label{eq:individual-balance}
 \Q(g_*)
  =\Q(f_*) - 2 c_*\Q(f_*,h) + c_*^2 \Q(h) = \Q(f_*) +  c_*^2\bigl(2a\E h+\Q(h)\bigr)
 =\Q(f_*)-c_*^2 S  ,
\end{equation}
and the sum over  $* \in \{D , \mathrm{sh} \}$ gives the identity in \eqref{eq:balance-identity}.

To prove the strict bound on the coefficients, use
$(k^{-1})'=-k'/k^2$ and integrate by parts in
\eqref{eq:linear-speeds} with respect to $k^{-2}(\theta) d\theta$ used in \eqref{eq:mu}. This gives
\begin{align*}
 \E f_D & =\E[(h - 2k) \cos 2 \theta] = \E[(a-k)\cos2\theta], \\
 \E f_{\mathrm{sh}} & =\E[(h-2k) \sin 2 \theta] = \E[(a-k)\sin2\theta].
\end{align*}
Since $a>0$ and $k>0$, it follows that
\[
 \sqrt{(\E f_D)^2+(\E f_{\mathrm{sh}})^2}
 =\left|\E\bigl[(a-k)(\cos2\theta,\sin2\theta)\bigr]\right|
 \leq\E|a-k|<\E(a+k)=\E h.
\]
Dividing by $\E h > 0$ completes the proof.
\end{proof}

\begin{corollary}[Linear instability criterion]\label{cor:linear-instability}
If a stationary oval satisfies $a>0$ and $S<0$, then it must be unstable via a linear deformation: there exists a symmetric matrix $B_0$ and an initial normal velocity $g = \ip{B_0 x}{N}$ with $\E g = 0$ such that $\Q(g) < 0$. 
If the oval is centrally symmetric, then $g$ is even and the deformation remains centrally symmetric. 
\end{corollary}
\begin{proof}
By Lemma~\ref{lem:balance}, if $S < 0$ then at least one of $g_D,g_{\mathrm{sh}}$ has
negative quadratic form, and we set $g = g_*$ for the corresponding $* \in \{D , {\mathrm{sh}}\}$. Note that $g = \ip{B_* x}{N} - c_* h = \ip{(B_* - c_* \Id) x}{N}$ is of the asserted form, and is even when the oval is centrally symmetric. 
The construction in Appendix~\ref{app:linear} realizes an exactly volume-preserving linear deformation with the same symmetry which reduces the Gaussian perimeter. 
\end{proof}

\subsection{The unconditional class}

The preceding corollary does not by itself establish instability within
the unconditional class, as the negative quadratic form might be the one resulting from the shear deformation. 
We treat this symmetry class separately.

\begin{lemma}[Unconditional instability criterion]
\label{lem:unconditional-instability}
Let $\Sigma$ be an unconditional stationary oval. If $a>0$ and $S<0$, then the oval is unstable within the
unconditional class. When $k'$ has no zero in $(0,\pi/2)$ (as for the remaining oval in Proposition \ref{prop:BCJ}), the
destabilizing family can be chosen among diagonal linear images of $\Sigma$ which preserve Gaussian volume (``unconditional squeezing''). 
\end{lemma}

\begin{proof}
Unconditionality implies $k'(0)=k'(\pi/2)=0$ by reflecting across the coordinate axes. It also makes $g_D$
unconditional, while $f_{\mathrm{sh}}$ changes sign under either
coordinate reflection, so $c_{\mathrm{sh}}=0$.

First suppose that $k'$ has no zero in $(0,\pi/2)$; by Lemma \ref{lem:period}, this is the case for alternative (iii) in Proposition \ref{prop:BCJ}. Choose
$w=\pm k'>0$ on this interval. Differentiating \eqref{eq:curvature-ode}
gives
\begin{equation}\label{eq:rotation-equation}
 w''+(1+k^{-2})w=0,
 \qquad w(0)=w(\pi/2)=0.
\end{equation}
The endpoint zeros are simple, since a multiplicity would force the solution to be constant zero by ODE uniqueness, a contradiction to noncircularity.
The function
$f_{\mathrm{sh}}$ vanishes at both endpoints by
\eqref{eq:linear-speeds}. Thus $f_{\mathrm{sh}}/w$ extends smoothly to
the closed interval. Writing $f_{\mathrm{sh}}=wz$ and using
\eqref{eq:rotation-equation},
\[
 (f_{\mathrm{sh}}')^2-(1+k^{-2})f_{\mathrm{sh}}^2
 =w^2(z')^2+(ww'z^2)'.
\]
The integral of the second term is zero (as $w$ vanishes at the endpoints), and hence
\begin{equation}\label{eq:ground-state}
 \int_0^{\pi/2}
 \bigl((f_{\mathrm{sh}}')^2-(1+k^{-2})f_{\mathrm{sh}}^2\bigr)\dd\theta
 =\int_0^{\pi/2}w^2
 \left(\left(\frac{f_{\mathrm{sh}}}{w}\right)'\right)^2\dd\theta
 \geq0.
\end{equation}
Reflection gives $\Q(f_{\mathrm{sh}})\geq0$ on the full curve. Therefore
\eqref{eq:balance-identity} implies
\[
 \Q(g_D)=(1-c_D^2)S-\Q(f_{\mathrm{sh}})<0.
\]
Since $g_D = f_D - c_D h = \ip{(B_D - c_D \Id) x}{N}$, it is the normal component of a diagonal linear map. The construction in Appendix~\ref{app:linear} realizes an exactly volume-preserving diagonal linear deformation preserving unconditionality and reducing Gaussian perimeter. 

Now suppose that $k'$ has a zero in $(0,\pi/2)$. This scenario is irrelevant for the proof of Theorem \ref{thm:main}, but is required for the proof of the very last sentence of Theorem \ref{thm:main2} about instability in the unconditional class. We sketch the proof, modeling it after \cite[Corollary~11.9]{Heilman} as in Remark \ref{rem:Heilman}. In this case, there are at least
two nodal intervals of $k'$ in this quadrant. Choose two of them and,
for $i=1,2$, define $u_i$ to equal $|k'|$ on the $i$th interval and zero
elsewhere in $[0,\pi/2]$. Extend each $u_i$ to the whole curve by using
coordinate reflections. These functions are nonzero, nonnegative,
unconditional and belong to $H^1(\Sigma)$. Note that by differentiating (\ref{eq:curvature-ode}), we have
\[
\LJac k' = k^2 k''' + (1 + k^2) k' = 0 . 
\]
Consequently, integration by parts on each supporting arc gives $\Q(u_i)=0$, since $k'$ vanishes at the arc endpoints. The interiors of their supports are disjoint, so
\[
 g:=(\E u_2)u_1-(\E u_1)u_2
 \quad\text{satisfies}\quad
 g\not\equiv0,\qquad \E g=0,\qquad \Q(g)=0.
\]
If the oval were stable in the unconditional class, smooth approximation
and polarization would give $\Q(g,\varphi)=0$ for every smooth
unconditional $\varphi$ with $\E\varphi=0$. Averaging test functions
over the coordinate reflections then yields $\LJac g \equiv C$ weakly for
some constant $C$. On the interior of a supporting arc, $g$ is a
constant multiple of $k'$, so $C=0$. The weak equation makes $g$ smooth.
But $g(0)=0$ by construction, and $g'(0)=0$ by reflection symmetry.
Uniqueness for the ODE $L_{Jac} g = 0$ now forces $g\equiv0$, a contradiction.
Thus there is a smooth unconditional mean-zero velocity with negative
quadratic form, and the construction in Appendix~\ref{app:stability}
realizes an exactly volume-preserving unconditional deformation which reduces the Gaussian perimeter. 
Note that in this scenario we do not assert that the deformation is linear.
\end{proof}

\section{The Key Inequality}
\label{sec:key-proof}

To apply the instability criteria from the previous section, it remains
to show that $a>0$ and $S<0$ for every noncircular stationary oval.
By \eqref{eq:S} and \eqref{eq:mu},
\[
 S=\frac{4 \pi + a\displaystyle\int_0^{2\pi} k^{-1}\dd\theta - 
 a^2\displaystyle\int_0^{2\pi} k^{-2}\dd\theta}
 {\displaystyle\int_0^{2\pi} k^{-2}\dd\theta}.
\]
Thus, to exclude any (noncircular) stationary oval as a stable candidate, it suffices to prove the following analytic inequality, which does not assume symmetry, stability or minimality of the corresponding
curve.

\begin{proposition}[Key Inequality]\label{prop:key}
Let $a\in\R$, and let $k:\R\to(0,\infty)$ be a nonconstant smooth
$2\pi$-periodic solution of
\begin{equation} \label{eq:k-ode}
 k''+k+a= \frac1k.
\end{equation}
Then $a>0$ and
\begin{equation}\label{eq:key-inequality}
 a^2\int_0^{2\pi}\frac{\dd\theta}{k^2}
 >a\int_0^{2\pi}\frac{\dd\theta}{k}+4\pi.
\end{equation}
\end{proposition}

\medskip

For the proof, it will be useful to introduce the following potential:
\begin{equation} \label{eq:V}
V(\kappa)=\frac{\kappa^2}{2}+a \kappa-\log \kappa,\qquad \kappa>0.
\end{equation}
Note that $V$ is strictly convex, and since $V(\kappa)$ tends to infinity at both endpoints of $(0,\infty)$, it must have a unique minimum at $b \in (0,\infty)$ given by 
\begin{equation} \label{eq:b}
0 = V'(b) = b + a - b^{-1} .
\end{equation}
The equation (\ref{eq:k-ode}) is precisely $k''+V'(k)=0$, and so multiplying by $k'$ gives
\[
\frac{d}{d\theta}
\left(\frac12(k')^2+V(k)\right)=0.
\]
Thus there is a preserved energy level $E \in \R$ such that
\begin{equation} \label{eq:energy-id}
\frac12(k')^2+V(k) \equiv E.
\end{equation}
Since $k$ is nonconstant, $E>V(b)$, as equality would force
both $k'=0$ and $V(k) \equiv V(b)$ everywhere. 

Strict convexity and the endpoint behavior of $V$ imply
that $V(\kappa)=E$ has exactly two roots $\kappa_-<b<\kappa_+$  in $(0,\infty)$, and that $V'(\kappa_-)<0<V'(\kappa_+)$.
The energy equation (\ref{eq:energy-id}) therefore gives
\[
\kappa_-\le k\le \kappa_+,\qquad
(k')^2=2(E-V(k)),
\]
and consequently
\[
k'=0
\quad\Longleftrightarrow\quad
k\in\{\kappa_-,\kappa_+\}.
\]
At these two turning values the curvature equation gives,
\[
k''=-V'(\kappa_-)>0,\qquad k''=-V'(\kappa_+)<0.
\]
Thus $\kappa_-$ can occur only at a strict minimum, and $\kappa_+$
only at a strict maximum. It follows that $k$ is strictly monotone
between successive extrema. Uniqueness for the ODE makes one excursion
from the maximum to the minimum and back a primitive cycle, and
time-reversal invariance gives equal lengths to its two branches.

\medskip

We first record the information about the primitive period needed for the Key Inequality and for Theorem~\ref{thm:main2}(i).

\begin{lemma}[Primitive period and parameter range]\label{lem:period-parameter}
Under the assumptions of Proposition~\ref{prop:key}, the primitive period
of $k$ is $2\pi/m$ for an integer $m\geq2$. Moreover, $a>0$, and
\begin{equation}\label{eq:period-equivalence}
 m=2\quad\Longleftrightarrow\quad 0<a<\frac2{\sqrt3}.
\end{equation}
\end{lemma}

\begin{proof}
Since $k$ is nonconstant and $2\pi$-periodic, its primitive period is
$2\pi/m$ for an integer $m\geq1$. By the preceding energy argument,
translating $\theta$ if necessary, we may assume that
\begin{equation}\label{eq:key-phase}
 k'(0)=k'(\pi/m)=0,
 \qquad k'<0\quad\text{on }I=(0,\pi/m).
\end{equation}
Differentiating \eqref{eq:k-ode}, and then differentiating $kk'$ twice,
gives
\begin{align}
 (k')''+(1+k^{-2})k'&=0,\label{eq:kprime-ode}\\
 (kk')''+4kk'&=(2/k-3a)k'.\label{eq:kkprime-ode}
\end{align}
Test both equations against $\sin(m\theta)$ on $I$. Integrating by parts
twice, the boundary terms vanish because $k'$, $kk'$ and
$\sin(m\theta)$ vanish at both endpoints. Hence the positive weight
$w=-k'\sin(m\theta)$ satisfies
\begin{equation}\label{eq:weighted-moments}
 \begin{aligned}
 \int_I\frac{w}{k^2}\dd\theta
 &=(m^2-1)\int_Iw\dd\theta,\\
 3a\int_Iw\dd\theta
 &=\int_I\left(\frac2k+(m^2-4)k\right)w\dd\theta.
 \end{aligned}
\end{equation}
The first identity excludes $m=1$, and the second then gives $a>0$.

If $m=2$, these identities become
\[
 \int_I\frac{w}{k^2}\dd\theta=3\int_Iw\dd\theta,
 \qquad
 \int_I\frac{w}{k}\dd\theta=\frac{3a}{2}\int_Iw\dd\theta.
\]
Since $1/k$ is nonconstant on $I$, Cauchy--Schwarz is strict:
\[
 \left(\frac{3a}{2}\right)^2
 =\frac{\left(\displaystyle\int_Iw/k\dd\theta\right)^2}
 {\left(\displaystyle\int_Iw\dd\theta\right)^2}
 <\frac{\displaystyle\int_Iw/k^2\dd\theta}
 {\displaystyle\int_Iw\dd\theta}=3.
\]
Thus $0<a<2/\sqrt3$.

Conversely, if $m\geq3$, then by AM-GM
\[
 \frac2k+(m^2-4)k\geq2\sqrt{2(m^2-4)}.
\]
Equality cannot hold throughout $I$, since $k$ is nonconstant.
Integrating against $w>0$ and using \eqref{eq:weighted-moments} gives
\begin{equation}\label{eq:higher-frequency-bound}
 a>\frac23\sqrt{2(m^2-4)}
 \geq\frac{2\sqrt{10}}3>\frac2{\sqrt3}.
\end{equation}
This proves \eqref{eq:period-equivalence}.
\end{proof}

We are now ready to prove the Key Inequality. All integrals below are with respect to Lebesgue measure in the normal
angle $\theta$.

\begin{proof}[Proof of Proposition \ref{prop:key}]
By Lemma~\ref{lem:period-parameter}, $a>0$. Set
\[
 F(\kappa):=\frac{a^2}{\kappa^2}-\frac a\kappa-2
 =\left(\frac a{\kappa^2}-\frac2\kappa\right)(\kappa+a).
\]
Our goal is to prove $\int_0^{2\pi}F(k)\dd\theta>0$. By Lemmas~\ref{lem:period} and~\ref{lem:period-parameter}, the reader only interested in the proof of Theorem \ref{thm:main} may skip Step 1 below, since the remaining oval in Proposition \ref{prop:BCJ} has primitive period $\pi$, and so necessarily $a < 2 / \sqrt3$. 
\medskip
\noindent\emph{Step 1. The range $a\geq2/\sqrt3$.}
Put $M=(2\pi)^{-1}\int_0^{2\pi}k^{-1}\dd\theta$.
Integrating \eqref{eq:k-ode} and using periodicity gives
$(2\pi)^{-1}\int_0^{2\pi}k\dd\theta=M-a$. Multiplying these two identities and applying Cauchy--Schwarz yields
\[
 M(M-a)\geq1,
  \text{ and hence }
 M\geq\frac{a+\sqrt{a^2+4}}2.
\]
If $a\geq2/\sqrt3$, this implies $aM\geq2$. Since $k^{-1}$ is
nonconstant, strict Cauchy--Schwarz gives
\[
 \frac1{2\pi}\int_0^{2\pi}F(k)\dd\theta
 >a^2M^2-aM-2=(aM-2)(aM+1)\geq0.
\]
This proves the Key Inequality in this range. We henceforth assume that
\begin{equation}\label{eq:a-bound}
 0<a<\frac2{\sqrt3},
\end{equation}
so by \eqref{eq:b},
\begin{equation}\label{eq:b-bound}
 \frac1{\sqrt3}<b=\frac{\sqrt{a^2+4}-a}{2}<1.
\end{equation}

\medskip
\noindent\emph{Step 2. A weighted inequality via Poincar\'e.}
We first prove the weighted inequality
$\int_0^{2\pi}F(k)(k')^2\dd\theta\geq0$. Introduce the function
\begin{equation}\label{eq:F-y}
 y:=(k+a)k'=\frac12\bigl((k+a)^2\bigr)'.
\end{equation}
In the range \eqref{eq:a-bound}, Lemma~\ref{lem:period-parameter}
shows that $k$ has primitive period $\pi$. Thus $y$ is $\pi$-periodic
and has zero mean on each period with respect to $\dd\theta$ (as a periodic derivative). The sharp Poincar\'e (or Wirtinger's) inequality on
the circle of length $\pi$ \cite[Theorem~258]{HardyLittlewoodPolya} (not $2 \pi$!)
therefore gives
\[
 \int_0^{2\pi}(y')^2\dd\theta = 2 \int_0^\pi (y')^2\dd\theta \geq 2 \cdot 4 \int_0^{\pi}y^2\dd\theta = 4 \int_0^{2\pi}y^2\dd\theta.
\]
Combining \eqref{eq:kprime-ode} and \eqref{eq:kkprime-ode} gives
\begin{equation}\label{eq:y-ode}
 y''+4y=\left(\frac2k-\frac a{k^2}\right)k'.
\end{equation}
Consequently, periodic integration by parts yields
\begin{equation}\label{eq:weighted-positive}
 \int_0^{2\pi}F(k)(k')^2\dd\theta
 =-\int_0^{2\pi}y(y''+4y)\dd\theta
 =\int_0^{2\pi}\bigl((y')^2-4y^2\bigr)\dd\theta\geq0.
\end{equation}

\medskip
\noindent\emph{Step 3. Removing the weight by using the energy identity.}
For $\kappa>0$ and $\kappa\ne b$, define
\begin{equation}\label{eq:R}
 R(\kappa)=\frac{F(\kappa)}{V'(\kappa)}
 =-\frac{2(\kappa+a)(\kappa-a/2)}
 {\kappa(\kappa-b)(\kappa+b^{-1})}.
\end{equation}
By \eqref{eq:b} and \eqref{eq:a-bound}, its poles and zeros alternate:
$-1/b<-a<0<a/2<b$.  It follows that all three residues are negative, namely (using $a = 1/b - b$):
\[
 R(\kappa)
 =-\frac{a^2}{\kappa}
 -\frac{3b^2-1}{b^2(1+b^2)}\,\frac1{\kappa-b}
 -\frac{b^2(3-b^2)}{1+b^2}\,\frac1{\kappa+1/b} ,
\]
where all three coefficients in front of the minus signs are positive by (\ref{eq:a-bound}) and (\ref{eq:b-bound}). It follows that
\begin{equation}\label{eq:Rprime}
 R'(\kappa)>0\qquad(\kappa>0,\ \kappa\ne b).
\end{equation}

At $b$, the function $V-V(b)$ has a double zero and $V'$ a simple zero.
Thus $R$ has at most a simple pole and $R'$ at most a double pole, so both
\[
 (V-V(b))R\qquad\text{and}\qquad(V-V(b))R'
\]
extend smoothly across $b$. The second product is nonnegative everywhere
and strictly positive away from $b$, by \eqref{eq:Rprime} and strict
convexity of $V$.

Using $k''=-V'(k)$ and $V'R=F$, differentiate the smooth periodic function
$(V(k)-V(b))R(k)k'$:
\begin{align*}
 &\bigl[(V(k)-V(b))R(k)k'\bigr]'\\
 &\qquad=\bigl[F(k)+(V(k)-V(b))R'(k)\bigr](k')^2
 -(V(k)-V(b))F(k).
\end{align*}
Substitute $V(k)-V(b)=E-V(b)-\tfrac12(k')^2$, which is
\eqref{eq:energy-id}, and integrate over $[0,2\pi]$. The integral of the
derivative on the left-hand side vanishes by periodicity, giving
\begin{equation}\label{eq:unweight}
 \begin{aligned}
 (E-V(b))\int_0^{2\pi}F(k)\dd\theta
 &=\frac32\int_0^{2\pi}F(k)(k')^2\dd\theta\\
 &\quad+\int_0^{2\pi}(V(k)-V(b))R'(k)(k')^2\dd\theta.
 \end{aligned}
\end{equation}
The first term on the right is nonnegative by \eqref{eq:weighted-positive}.
The second is strictly positive: its integrand is nonnegative everywhere
and positive wherever $k'\ne0$ and $k\ne b$. These conditions fail only
at the two turning points and the two crossings of $b$ in each primitive
period. Since $E-V(b)>0$, we conclude that
$\int_0^{2\pi}F(k)\dd\theta>0$, establishing \eqref{eq:key-inequality}.
\end{proof}

Consequently every stationary oval satisfies $a>0$ and
$S<0$. In the next section, we combine this with the
instability criteria of Section~\ref{sec:stability} to prove
Theorems~\ref{thm:main} and~\ref{thm:main2}.

\section{Proofs of Theorems \ref{thm:main} and \ref{thm:main2}}
\label{sec:completion}

\begin{proof}[Proof of Theorem~\ref{thm:main}]
Fix $v\in(0,1)$ and take an unconditional minimizer. By
Proposition~\ref{prop:BCJ}, after possible complementation, it must be a symmetric strip (aligned with one of the coordinate axes),
a centered ball, or a stationary oval; by Lemma~\ref{lem:period}, in the last case it must have primitive period $\pi$. 
The Key Inequality of Proposition~\ref{prop:key}, together with the
unconditional instability criterion of
Lemma~\ref{lem:unconditional-instability}, excludes the oval (showing that an appropriate linear ``unconditional squeezing'' will preserve Gaussian volume while strictly reducing Gaussian perimeter). 
Thus every unconditional minimizer is a symmetric coordinate strip, a centered disk, or a complement of one of these sets. 

Lemma~\ref{lem:unconditional-2} proves equality of the centrally symmetric and
unconditional Gaussian isoperimetric problems in $\R^2$:
\[
\inf P_\gamma(E) = \inf  P_\gamma(F) ,
\]
where the first infimum is over all centrally symmetric measurable sets $E$, and the second is over all unconditional measurable sets $F$, satisfying $\gamma(E) = \gamma(F) = v$ for a prescribed $v \in (0,1)$. 
It remains to classify all centrally symmetric minimizers, rather than only to produce an unconditional one.

Let $E=-E$ be an arbitrary minimizer and put $p=P_\gamma(E)$.
Choose equipartitioning perpendicular axes as in
Lemma~\ref{lem:unconditional-2}, and use the coordinates determined by
these axes. Denote the four open quadrants by
$Q_1,\ldots,Q_4$ in cyclic order. For each $j$, reflect $E\cap Q_j$
into all four quadrants, obtaining $F_j$. Since $\gamma(F_j) = 4 \gamma(E \cap Q_j) = \gamma(E)$, minimality and
\eqref{eq:reflection-exact} yield
\[
 p\le P_\gamma(F_j)=4P_\gamma(E;Q_j).
\]
Summing over $j$ gives
\[
 4p\le4\sum_{j=1}^4P_\gamma(E;Q_j)\le 4 P_{\gamma}(E) = 4 p.
\]
Hence equality holds throughout, and so every $F_j$ is an unconditional
minimizer, and $E$ has no perimeter on either cutting axis. The BV jump
formula therefore implies agreement of the one-sided traces of $E$
across each open half-axis.

Central symmetry gives $F_1=F_3$ and $F_2=F_4$; denote these two unconditional minimizing
sets by $A$ and $B$. Agreement of the traces across the two positive
half-axes implies that $A$ and $B$ have the same ordered pair of traces
there. We now use that such a pair determines the candidate uniquely: indeed, writing
the positive half-axis as $(0,\infty)$, our initial classification of unconditional minimizers shows that the only  possibilities are
\[
\begin{array}{c|c}
 \text{candidate}&\text{traces on the two positive half-axes}\\ \hline
 B_r^2&((0,r),(0,r))\\
 \R^2\setminus B_r^2&((r,\infty),(r,\infty))\\
 \{|x_1|<r\}&((0,r),(0,\infty))\\
 \{|x_2|<r\}&((0,\infty),(0,r))\\
 \{|x_1|>r\}&((r,\infty),\varnothing)\\
 \{|x_2|>r\}&(\varnothing,(r,\infty)).
\end{array}
\]
Note that endpoints are irrelevant for traces. Since $v \in (0,1)$, all radii are
finite and positive, and so distinct entries cannot give the same ordered
trace pair. We conclude that $A=B$. On every quadrant $E$ agrees with this same
candidate, so $E=A$ modulo a null set. Rotating the equipartitioning axes
back proves the full classification.
\end{proof}

\begin{proof}[Proof of Theorem~\ref{thm:main2}]
As explained in Section~\ref{sec:structure}, $\Sigma$ is strictly convex, and its curvature
$k$, parametrized by the outward normal angle, is a nonconstant positive
$2\pi$-periodic solution of the curvature equation (\ref{eq:curvature-ode}). Part (i) therefore
follows from Lemma~\ref{lem:period-parameter}.

For part (ii), Proposition~\ref{prop:key} gives $a>0$ and, by
(\ref{eq:S}) and (\ref{eq:mu}), $S<0$. Corollary~\ref{cor:linear-instability}
therefore supplies a symmetric matrix $B_0$ such that
$g=\ip{B_0x}{N}$ satisfies $\E g=0$ and $\Q(g)<0$, proving instability (in the centrally symmetric class if $\Sigma$ is centrally symmetric). As explained in the appendix, this velocity is realized by an exactly Gaussian-volume-preserving linear deformation. 

If $\Sigma$ is unconditional, the stronger conclusion of instability within the
unconditional class follows from Lemma~\ref{lem:unconditional-instability}. As follows from its proof, when the primitive period is $\pi$, the resulting velocity is of the form $g=\ip{B_0x}{N}$ for a diagonal $B_0$, and again it is realized by an exactly Gaussian-volume-preserving linear diagonal deformation. However, when the primitive period is strictly smaller than $\pi$, we make no such claim. 
\end{proof}

\appendix
\section{Second variation and an exact volume-preserving deformation}
\label{app:variation}

In this appendix, which is provided for completeness, we compute the standard variation formulas with the normalization used in
this work, and then construct the required volume-preserving family of competitors explicitly.
Let $K\subset\R^2$ be a smooth bounded domain with Gaussian stationary boundary
$\Sigma$, namely satisfying $k-h \equiv \lambda \in \R$ with $h = \ip{x}{N}$.

\subsection{Computation of second variation along normal deformation}

For a smooth normal velocity $f$, consider the map $\Sigma \ni x\mapsto x+tfN \in \R^2$,
with the original curve parametrized by arclength. Since
$\partial_sN=kT$, its line element and Gaussian density satisfy
\begin{align*}
 \frac{\dd s_t}{\dd s}
 &=1+tkf+\frac{t^2}{2}|\partial_sf|^2+O(t^3),\\
 \frac{\rho_2(x+tfN)}{\rho_2(x)}
 &=1-thf+\frac{t^2}{2}(h^2-1)f^2+O(t^3).
\end{align*}
The tubular-coordinate volume element is
$(1+rk)\rho_2(x+rN)\dd r\dd s$. Therefore
\begin{align}
 \gamma'(0)&=\int_\Sigma f\rho_2\dd s,
 &P_\gamma'(0)&=\int_\Sigma(k-h)f\rho_2\dd s,
 \label{eq:first-variation}\\
 \gamma''(0)&=\lambda\int_\Sigma f^2\rho_2\dd s,
 &P_\gamma''(0)&=\int_\Sigma
 \bigl(|\partial_sf|^2+(h^2-1-2kh)f^2\bigr)\rho_2\dd s.
 \label{eq:second-normal}
\end{align}
Subtracting $\lambda\gamma''(0)$ and using $(k-h)^2=\lambda^2$ yields
\begin{equation}\label{eq:second-lagrangian}
  (P_\gamma-\lambda\gamma)''(0)
 =\int_\Sigma\bigl(|\partial_sf|^2-(1+k^2)f^2\bigr)\rho_2\dd s
 =-\int_\Sigma f\LJac f\rho_2\dd s = P_\gamma(K) \Q(f) . 
\end{equation}
The drift term in the definition \eqref{eq:Jacobi-intrinsic} of $L_{Jac}$ is exactly the term supplied
by weighted integration by parts.

The Lagrangian $P_\gamma-\lambda\gamma$ has zero first variation at
$\Sigma$. Consequently its second derivative along any smooth family
depends only on the initial normal speed: an acceleration term is paired
with the vanishing first variation, and tangential motion is a
reparametrization. 
Now, let $K_t$ be any smooth family with $K_0=K$ and
initial normal velocity $f$. By stationarity of the Lagrangian, the
normal-graph computation \eqref{eq:second-lagrangian} gives
\[
 \left.\frac{d^2}{dt^2}P_\gamma(K_t)\right|_{t=0}
 -\lambda\left.\frac{d^2}{dt^2}\gamma(K_t)\right|_{t=0}
 =P_\gamma(K)\Q(f).
\]
If the family preserves Gaussian volume exactly, namely
$\gamma(K_t)=\gamma(K)$ for all sufficiently small $t$, the second
term on the left vanishes. Hence
\[
 \left.\frac{d^2}{dt^2}P_\gamma(K_t)\right|_{t=0}
 =P_\gamma(K)\Q(f).
\]

\subsection{Justification of the stability implication}
\label{app:stability}

Assume that $K$ minimizes Gaussian perimeter among 
sets of Gaussian volume $\gamma(K)$ in a certain symmetry class (e.g. centrally symmetric or unconditional).
Let $f\in C^\infty(\Sigma)$ belong to the same symmetry class and
satisfy $\E f=0$. We construct an exactly volume-preserving family
with initial normal velocity $f$ within the same symmetry class.

For $t,r$ sufficiently small, let $K_{t,r}$ be the normal perturbation
of $K$ with boundary
\[
 \partial K_{t,r}
 =\{x+(tf(x)+r)N(x):x\in\Sigma\}.
\]
These domains are smooth and in the same symmetry class as $K$ (by invariance of $f$ and equivariance of $N$). Define:
\[
 \Psi(t,r)=\gamma(K_{t,r})-\gamma(K).
\]
By the first-variation formula \eqref{eq:first-variation},
\[
 \partial_r\Psi(0,0)=\int_\Sigma\rho_2\dd s=P_\gamma(K)>0,
 \qquad
 \partial_t\Psi(0,0)=\int_\Sigma f\rho_2\dd s
 =P_\gamma(K)\E f=0.
\]
The implicit function theorem therefore gives a smooth function
$\alpha(t)$ such that
\[
 \alpha(0)=0,\qquad
 \Psi(t,\alpha(t))=0,\qquad
 \alpha'(0)=-\frac{\partial_t\Psi(0,0)}
 {\partial_r\Psi(0,0)}=0.
\]
Consequently, $K_t:=K_{t,\alpha(t)}$ has constant Gaussian volume,
and its initial normal velocity is $f$, because $\alpha'(0)=0$.

Minimality gives $P_\gamma(K_t)\ge P_\gamma(K)$ for all sufficiently
small $t$. Since $\gamma(K_t)$ is constant, the second-variation
identity established above yields
\[
 0\le
 \left.\frac{d^2}{dt^2}P_\gamma(K_t)\right|_{t=0}
 =P_\gamma(K)\Q(f).
\]
As $P_\gamma(K)>0$, this proves \eqref{eq:stability}.

\subsection{A linear family of deformations} \label{app:linear}

When $f(x) = \ip{B x}{N}$, a simple modification of the above construction yields a linear family of deformations. Indeed, simply use
\[
K_{t,r} = e^{r} e^{t B} K ,
\]
and repeat the above argument. Now $\partial_r\Psi(0,0)=P_\gamma(K)\E h>0$, since
$h=k+a>0$ for the stationary ovals considered above.  The resulting  family $K_t:=K_{t,\alpha(t)} = e^{\alpha(t)} e^{t B} K$ has constant Gaussian volume,
and its initial normal velocity is $f$.  If $K$ is centrally symmetric then so are $K_{t}$, and when $K$ is unconditional and $B$ is diagonal, then $K_{t}$ are unconditional as well.


\begin{thebibliography}{99}

\bibitem{AbreschLanger}
U.~Abresch and J.~Langer,
\emph{The normalized curve shortening flow and homothetic solutions},
J. Differential Geom. \textbf{23} (1986), no.~2, 175--196.

\bibitem{AFP}
L.~Ambrosio, N.~Fusco and D.~Pallara,
\emph{Functions of Bounded Variation and Free Discontinuity Problems},
Oxford Mathematical Monographs, Oxford University Press, Oxford, 2000.

\bibitem{BCJ}
M.~Barchiesi, F.~Cagnetti and V.~Julin,
\emph{Double-arc rearrangement and the symmetric isoperimetric problem
in Gauss space}, preprint, 2026,
\href{https://arxiv.org/abs/2610.15133v1}{arXiv:2610.15133v1}.

\bibitem{BarchiesiJulin}
M.~Barchiesi and V.~Julin,
\emph{Symmetry of minimizers of a Gaussian isoperimetric problem},
Probab. Theory Related Fields \textbf{177} (2020), 217--256.

\bibitem{Barthe}
F.~Barthe,
\emph{An isoperimetric result for the Gaussian measure and unconditional sets},
Bull. London Math. Soc. \textbf{33} (2001), no.~4, 408--416.

\bibitem{Borell}
C.~Borell,
\emph{The Brunn--Minkowski inequality in Gauss space},
Invent. Math. \textbf{30} (1975), no.~2, 207--216.

\bibitem{ChakrabartiRegev}
A.~Chakrabarti and O.~Regev,
\emph{An optimal lower bound on the communication complexity of
Gap-Hamming-Distance},
SIAM J. Comput. \textbf{41} (2012), no.~5, 1299--1317.

\bibitem{Chang}
J.-E.~Chang,
\emph{1-Dimensional solutions of the $\lambda$-self shrinkers},
Geom. Dedicata \textbf{189} (2017), 97--112.

\bibitem{ChengWei}
Q.-M.~Cheng and G.~Wei,
\emph{Complete $\lambda$-hypersurfaces of weighted volume-preserving mean curvature flow},
Calc. Var. Partial Differential Equations \textbf{57} (2018), no.~2,
Paper No.~32, 21 pp.

\bibitem{FuscoLaManna}
N.~Fusco and D.~A.~La Manna,
\emph{A remark on a conjecture on the symmetric Gaussian problem},
Proc. Edinb. Math. Soc. (2) \textbf{67} (2024), no.~2, 643--661.

\bibitem{HardyLittlewoodPolya}
G.~H.~Hardy, J.~E.~Littlewood and G.~P\'olya,
\emph{Inequalities}, 2nd ed.,
Cambridge University Press, Cambridge, 1952.

\bibitem{Heilman}
S.~Heilman,
\emph{Symmetric convex sets with minimal Gaussian surface area},
Amer. J. Math. \textbf{143} (2021), no.~1, 53--94. 

\bibitem{LaManna}
D.~A.~La Manna,
\emph{Local minimality of the ball for the Gaussian perimeter},
Adv. Calc. Var. \textbf{12} (2019), no.~2, 193--210.

\bibitem{LatalaOleszkiewicz}
R.~Lata{\l}a and K.~Oleszkiewicz,
\emph{Gaussian measures of dilatations of convex symmetric sets},
Ann. Probab. \textbf{27} (1999), no.~4, 1922--1938.

\bibitem{Maggi}
F.~Maggi,
\emph{Sets of Finite Perimeter and Geometric Variational Problems:
An Introduction to Geometric Measure Theory},
Cambridge Studies in Advanced Mathematics, vol.~135,
Cambridge University Press, Cambridge, 2012.

\bibitem{McGonagleRoss}
M.~McGonagle and J.~Ross,
\emph{The hyperplane is the only stable, smooth solution to the
isoperimetric problem in Gaussian space},
Geom. Dedicata \textbf{178} (2015), 277--296.

\bibitem{EMilmanCubeCPAM}
E.~Milman,
\emph{Isoperimetric inequalities on slabs with applications to cubes and
Gaussian slabs},
Comm. Pure Appl. Math. \textbf{79} (2026), no.~4, 1012--1072.

\bibitem{EMilmanNeemanAnnals}
E.~Milman and J.~Neeman,
\emph{The Gaussian double-bubble and multi-bubble conjectures},
Ann. of Math. (2) \textbf{195} (2022), no.~1, 89--206.

\bibitem{Musin}
O.~R.~Musin,
\emph{Borsuk–Ulam type theorem for the orthogonal group and orthogonal four-partitions},
preprint, 2026,
\href{https://arxiv.org/abs/2609.09688v3}{arXiv:2609.09688}.

\bibitem{ODonnell}
R.~O'Donnell,
\emph{Open problems in analysis of Boolean functions},
preprint, 2012,
\href{https://arxiv.org/abs/1204.6447}{arXiv:1204.6447}.

\bibitem{SudakovTsirelson}
V.~N.~Sudakov and B.~S.~Tsirel'son,
\emph{Extremal properties of half-spaces for spherically invariant measures},
Zap. Nau\v{c}n. Sem. Leningrad. Otdel. Mat. Inst. Steklov. (LOMI)
\textbf{41} (1974), 14--24, 165 (Russian).


\end{thebibliography}
\end{document}